\documentclass{elsart3-1}
\usepackage{mathrsfs} 
\usepackage{amsmath} 

\usepackage{amscd,amssymb,graphics,lineno}
\usepackage{mathrsfs} 
\usepackage{amsopn}

\usepackage{amssymb}

\usepackage[english,french]{babel}

\newtheorem{theorem}{Theorem}
\newtheorem{lemma}[theorem]{Lemma}
\newtheorem{e-proposition}[theorem]{Proposition}
\newtheorem{corollary}[theorem]{Corollary}
\newtheorem{e-definition}[theorem]{Definition\rm}

\newtheorem{question}{\it Question\/}
\newtheorem{theoreme}[theorem]{Th\'eor\`eme}

\newtheorem{proposition}[theorem]{Proposition}
\newtheorem{corollaire}[theorem]{Corollaire}

\renewcommand{\theequation}{\arabic{equation}}
\def\og{\leavevmode\raise.3ex\hbox{$\scriptscriptstyle\langle\!\langle$~}}
\def\fg{\leavevmode\raise.3ex\hbox{~$\!\scriptscriptstyle\,\rangle\!\rangle$}}

\newcommand{\abs}[1]{\vert#1\vert}
\def\norm#1{\left\Vert#1\right\Vert}
\def\R {{\mathbb R}}

\def\I {{\mathbb I}}

\def\N{{\mathbb N}}

\def\e{{\varepsilon}}

\def\s{{\mathbb S}}

\def\H{{\mathcal H}}

\def\Aut{{\mathrm{Aut}\,}}
\def\Ad{{\mathrm{Ad}\,}}

\def\RUCB{{\mbox{\rm RUCB}\,}}
\def\LUCB{{\mbox{\rm LUCB}\,}}

\begin{document}

\centerline{Group theory/Th\'eorie des groupes} 

\begin{frontmatter}


\selectlanguage{english}
\title
{An amenability-like property of finite energy path and loop groups
\\[3mm]
\large\sl 
Une propri\'et\'e semblable \`a la moyennabilit\'e des groupes de chemins et de lacets \`a \'energie finie
}

\selectlanguage{english}
\author
{Vladimir Pestov}
\ead{vladimir.pestov@uottawa.ca}

\address
{Departamento de Matem\'atica,
Universidade Federal de Santa Catarina,
Campus Universit\'ario Trindade,
CEP 88.040-900 Florian\'opolis-SC, Brasil 
}
\address{
D\'epartement de math\'ematiques et de statistique,
Universit\'e d'Ottawa, Complexe STEM, 150 Louis-Pasteur Pvt,
Ottawa, Ontario K1N 6N5 Canada
}
\medskip

\begin{abstract}
\selectlanguage{english}
We show that the groups of finite energy loops and paths (that is, those of Sobolev class $H^1$) with values in a compact connected Lie group, as well as their central extensions, satisfy an amenability-like property: they admit a left-invariant mean on the space of bounded functions uniformly continuous with regard to a left-invariant metric.
Every strongly continuous unitary representation $\pi$ of such a group (which we call skew-amenable) has a conjugation-invariant state on $B({\mathcal H}_{\pi})$.

\vskip 0.5\baselineskip

\selectlanguage{french}
\noindent{\bf R\'esum\'e} \vskip 0.5\baselineskip \noindent
Nous montrons que les groupes de lacets et de chemins \`a \'energie finie (c.\`a.d. de classe $H^1$ de Sobolev) \`a valeurs dans un groupe de Lie compact et connexe, ainsi que leurs extensions centrales, satisfont une version de la moyennabilit\'e: ils admettent une moyenne invariante \`a gauche sur l'espace de fonctions born\'ees uniform\'ement continues par rapport a une m\'etrique invariante \`a gauche. Chaque repr\'esentation unitaire continue, $\pi$, d'un tel groupe (que nous disons d'\^etre ``moyennable en biais'') poss\`ede un \'etat sur $B({\mathcal H}_{\pi})$ invariant sous conjugaison.
\end{abstract}
\end{frontmatter}

\selectlanguage{french}
\section*{Version fran\c{c}aise abr\'eg\'ee}

\subsubsection*{Motivation et r\'esultats.}
On s'int\'eresse \`a la question suivante de Carey et Grundling (\cite{CG}, p. 114). Soit $X$ une vari\'et\'e riemannienne compacte. Les groupes $C^k(X,SU(n))$, $C^{\infty}(X,SU(n))$, et $H^k(X,SU(n))$, munis des topologies correspondantes, sont-ils moyennables? Cette question, \`a son tour, est orient\'ee vers le probl\`eme de l'existence d'un \'etat de plus basse \'energie d'une th\'eorie de jauge. 

On rappelle qu'un groupe topologique $G$ (pas forc\'ement localement compact) est dit moyennable (au sens de Pierre de la Harpe \cite{dlH2}) si toute action continue de $G$ sur un espace compact quelconque poss\`ede une mesure de probabilit\'e borelienne reguli\`ere et invariante. Cela \'equivaut \`a l'existence d'une moyenne invariante \`a gauche sur l'espace $\RUCB(G)$ des fonctions born\'ees uniform\'ement continues \`a droite. Ici, nous suivons Bourbaki (\cite{bourbaki}, Ch. III, p. 20) en disant qu'une fonction $f\colon G\to\R$ est uniform\'ement continue \`a droite (respectivement, \`a gauche) si pour chaque $\e>0$ il existe un voisinage $V$ de l'unit\'e de $G$ tel que, quels que soient $x,y\in G$ avec $xy^{-1}\in V$ (resp., $x^{-1}y\in V$), on a $\abs{f(x)-f(y)}<\e$. 

\`A notre connaissance, la question de Carey et Grundling reste ouverte pour toutes les valeurs de $k\geq 0$ (on exige $k> \dim X/2$ pour que les groupes de Sobolev soient bien d\'efinis). La seule exception est le cas $k=0$ dans la dimension un, c'est-\`a-dire, le cas des groupes de chemins et de lacets continus, dont la moyennabilit\'e est une cons\'equence imm\'ediate des r\'esultats de Marie-Paule Malliavin e Paul Malliavin \cite{MM} (voir ci-apr\`es la section \ref{s:mm}).

Dans cette Note, nous nous occupons des groupes de chemins et lacets d'une classe de r\'egularit\'e l\'eg\`erement plus \'elev\'ee, \`a savoir, de la classe de Sobolev $H^1$ (aussi connus comme des chemins et lacets \`a \'energie finie). Bien que nous ne sommes pas en mesure de montrer la moyennabilit\'e des groupes en question, nous \'etablissons une variation de cette propri\'et\'e \'etroitement li\'ee et, de toute \'evidence, bien pertinante.

Disons qu'un groupe topologique $G$ est {\em moyennable en biais}\footnote{Cette suggestion terminologique, visant l'expression anglaise {\em skew-amenable group,} appartient \`a Colas Bardavid.} s'il existe une moyenne invariante \`a gauche sur les fonctions born\'ees uniform\'ements continues \`a gauche. Th\'eor\`eme 2.2.1 dans \cite{greenleaf} entraine que pour les groupes localement compacts, la moyennabilit\'e en biais et la moyennabilit\'e sont deux notions \'equivalentes. Ceci est inexact sans l'hypoth\`ese de la compacit\'e locale: le groupe unitaire $U(\ell^2)$ muni de la topologie forte est moyennable \cite{dlH2}, mais pas moyennable en biais (exemple 3.6.3 dans \cite{P06}). En revanche, nous ne savons pas si tout groupe moyennable en biais est moyennable.

Voici le r\'esultat principal de la Note.

 \vskip .2cm
\begin{theoreme}
Soit $G$ un groupe de Lie compact et connexe. Les groupes de chemins et de lacets \`a \'energie finie \`a valeurs dans $G$, ainsi que leurs variantes bas\'ees et extensions centrales, sont moynnables en biais.
\label{th:principal}
\end{theoreme}
\vskip .2cm

Tout groupe topologique $G$ moyennable en biais poss\`ede la propr\'et\'e suivante. Soit $\pi$ une repr\'esentation unitaire et fortement continue de $G$ dans un espace hilbertien ${\mathcal H}_{\pi}$. Alors, l'alg\`ebre d'operateurs $B({\mathcal H}_{\pi})$ admet un \'etat $\phi$ invariant par conjugaisons: quels que soient $T\in B({\mathcal H}_{\pi})$ et $g\in G$, on a $\phi(\pi_{g}^{\ast} T\pi_g)=\phi(T)$. En d'outres mots, la repr\'esentation $\pi$ est moyennable au sens de Bekka \cite{Bek1}. Pour la preuve, voir \cite{GP07}, Prop. 4.5.
On en d\'eduit:
\vskip .2cm

\begin{corollaire}
Chaque repr\'esentation unitaire et continue, $\pi$, d'un groupe mentionn\'e dans le th\'eor\`eme \ref{th:principal} poss\`ede un \'etat invariant $\phi$ sur l'alg\`ebre $B({\mathcal H}_{\pi})$.
\end{corollaire}
\vskip .2cm

\`A notre avis, il serait profitable d'enqu\^eter sur la moyennabilit\'e en biais des groupes vus dans la question de Carey et Grundling au d\'ebut de la Note.

\subsubsection*{Esquisse de la preuve.}
Soit $G$ un groupe de Lie compact r\'eel d'alg\`ebre de Lie $\mathfrak g$. On r\'ealise $G$ et $\mathfrak g$ par des matrices, et munit $M_n(\R)$ de la norme de Hilbert--Schmidt, not\'ee $\norm{\cdot}_2$.
Toute application $f\colon [0,1]\to G$ absolument continue poss\`ede la {\em deriv\'ee logarithmique \`a droite} d\'efinie presque partout par $\partial^{log}f(t)=f^{\prime}(t)\cdot f(t)^{-1}$. Si $f$ v\'erifie la condition $\int_0^1\norm{\partial^{log}f(t)}\,dt<\infty$, on dit que $f$ est un chemin \`a \'energie finie. Cela \'equivaut \`a $f$ \^etre de la classe de Sobolev $H^1$. L'ensemble de toutes les chemins \`a \'energie finie \`a valeurs dans $G$, muni de la topologie $H^1$, forme un groupe topologique par rapport \`a la multiplication ponctuelle (en effet, un groupe de Lie banachique, voir \cite{AHKMT}, Sect. 1.8), not\'e $H^1(\I,G)$. De m\^eme, pour le groupe de lacets, $H^1(\s^1,G)$.

La restriction de la d\'eriv\'ee logarithmique \`a droite sur le sous-groupe $H^1_0(\I,G)$ des chemins bas\'es (c.\`a.d., ceux-ci v\'erifiant $f(0)=e$) \'etablit un hom\'eomorphisme entre le groupe en question et l'espace hilbertien $L^2(\I,{\mathfrak g})$, l'application r\'eciproque \'etant donn\'ee par l'int\'egrale multiplicative \cite{DF}, \begin{equation}
t\mapsto \prod_0^t\exp f(s)\,ds.
\label{eq:prodexp}
\end{equation}

La d\'eriv\'ee logarithmique \'etant un cocycle (voir l'\'equation (\ref{eq:cocycle})), elle \'etablit un isomorphisme entre le groupe topologique $H^1_0(\I,G)$ et l'espace $L^2(\I,{\mathfrak g})$ muni de la loi multiplicative
\[f\ast g = f + \Ad_{\Pi\exp f}g.\]
Ici, $\Pi\exp f$ d\'esigne la fonction dans l'\'equation (\ref{eq:prodexp}), et $\Ad$ est la repr\'esentation adjointe, appliqu\'ee ponctuellement (pour les groupes matriciels, il s'agit tout simplement de la conjugaison). On note, en particulier, que la structure uniforme additive de l'espace hilbertien $L^2(\I,{\mathfrak g})$ est la structure uniforme \`a gauche du groupe topologique $(L^2(\I,{\mathfrak g}),\ast)$, parce que elle est invariante par translations et rotations. Par cons\'equent, dans le but de montrer la moyennabilit\'e en biais du groupe $H^1_0(\I,G)$ --- o\`u, de m\^eme, de sa copie isomorphe $(L^2(\I,{\mathfrak g}),\ast)$ ---, il suffit de construire une moyenne invariante sur l'espace des fonctions born\'ees et uniform\'ement continues sur $L^2(\I,{\mathfrak g})$ invariante par toutes les translations, ainsi que par les rotations sous la forme sp\'eciale $\Ad_f$, o\`u $f\in H^1_0(\I,G)$.

Pour chaque entier positif, $N$, notons $V_N$ le sous-espace euclidien de $L^2(\I,{\mathfrak g})$ de dimension finie $N\cdot\dim{\mathfrak g}$ form\'e par des fonctions simples prenant une valeur constante sur chaque interval $[i/N,(i+1)N)$. Notons $B_N$ la boule unit\'e ferm\'ee dans l'espace $V_N$.
Soit $(R_N)$ une suite des r\'eels positifs convergeant vers l'infini. Notons $\nu_N$ la restriction de la mesure de Lebesgue sur la boule $R_NB_N$, normalis\'ee pour devenir une mesure de probabilit\'e. (L'on peut aussi bien utiliser une measure gaussienne centr\'ee de variance $R_N^2$). Fixons un ultrafiltre libre $\mathcal U$ sur $\N$, et d\'efinissons, pour chaque fonction $F$ born\'ee et uniform\'ement continue sur $L^2(\I,{\mathfrak g})$,
\[\phi(F)=\lim_{N\to{\mathcal U}}\int_{V_N} F(x)\,d\nu_N(x).\]
La fonctionnelle $\phi$ est lin\'eaire, positive, et de norme $1$. Elle est invariante par translations d\`es que $R_N=\omega(N^{1/2})$. Sous cette hypoth\`ese, un raisonnement g\'eom\'etrique simple montre que la valeur de la mesure de la diff\'erence sym\'etrique $B_N\Delta (f_N+B_N)$ converge vers z\'ero quand $N\to\infty$ pour toute suite $(f_N)$ born\'ee dans $L^2(\I,{\mathfrak g})$ et telle que $f_N\in V_N$. Maintenant on use la densit\'e de sous-espace vectoriel $\cup_NV_N$ dans $L^2(\I,{\mathfrak g})$.

La preuve de l'invariance par les transformations $\Ad_r$, o\`u $r\in H^1_0(\I,G)$, est l\'eg\`erement plus compliqu\'ee. Elle a lieu sous l'hypoth\`ese $R_N= o(N/\log N)$. \'Etant donn\'e une fonction $r$ ci-dessus, on l'approxime avec une fonction $\rho\colon\I\to G$ simple, prenant des valeurs constantes sur les intervalles $[i/N,(i+1)/N)$. La transformation $\rho$ conserve la mesure $\nu_N$, d'o\`u il suffit d'estimer la norme de $(\Ad_r-\Ad_{\rho})f$ pour un \'el\'ement quelconque $f$ de $R_NB_N$. 

En utilisant la continuit\'e $1/2$-h\"olderienne de $r$, on obtient
\[\norm{\Ad_rf-\Ad_{\rho} f}_2\leq CN^{-1/2}\norm{\partial^{\log}r}
\max_{i=0}^{N-1}\left(N^{-1/2}\norm{f_i}_2 \right),\]
o\`u $f_i\in {\mathfrak g}$ est la valeur constante de $f$ sur l'intervalle $i$-i\`eme. \`A l'aide de l'in\'egalit\'e de concentration de Paul L\'evy pour la boule euclidien (\cite{ledoux}, Prop. 2.9, p. 30), on montre que la mesure $\nu_N$ de l'ensemble de toutes les fonctions $f$ ayant la propri\'et\'e
$\max_{i=0}^{N-1}\left(N^{-1/2}\norm{f_i}_2 \right)<R_NN^{-1/2}\log N$ converge vers $1$ quand $N\to\infty$. D'o\`u on d\'eduit que 
\[\int_{L^2} \left[F(\Ad_rf) - F(f)\right] d\nu_N(f)\to 0\]
pour toute fonction $F$ born\'ee et uniform\'ement continue sur $L^2(\I,{\mathfrak g})$. On vient donc d'\'etablir que le groupe de chemins bas\'es $H^1_0(\I,G)$ est moyennable en biais.

Pour transf\'erer le r\'esultat au groupe de lacets bas\'es $H^1_0(\s^1,G)$, nous montrons que tout sous-groupe normal et co-compact d'un groupe moyennable en biais est moyennable en biais. Notons que en g\'en\'eral, la moyennabilit\'e en biais n'est pas transmise au sous-groupes topologiques ferm\'es, comme le t\'emoigne l'exemple dans \cite{CT}.

Le r\'esultat pour les groupes de chemins et de lacets est en d\'eduit avec l'aide du th\'eor\`eme suivant: le produit semi-direct d'un groupe compact agissant sur un groupe moyennable en biais est moyennable en biais. Le r\'esultat est donc plus faible et l'argument, plus d\'elicat que pour les groupes moyennables, d\^u au fait que l'action du groupe \`a gauche sur les fonctions uniform\'ement continues \`a gauche est g\'en\'eralement discontinue. Finalment, un argument beaucoup plus simple \'etablit le r\'esultat pour les extensions centrales.
Effectivement, nous ne savons pas si la classe des groupes moyennables en biais est ferm\'ee par extensions les plus g\'en\'erales.

\selectlanguage{english}
\section{Motivation and statements of results}
\label{}

We are interested in the following open question.

\vskip .2cm
\begin{question}[Carey and Grundling \cite{CG}, p. 114] 
Let $X$ be a compact Riemannian manifold. Are the groups $C^k(X,SU(n))$, $C^{\infty}(X,SU(n))$, and $H^k(X,SU(n))$ with their natural group topologies amenable?
\end{question}
\vskip .2cm

Amenability of a topological group $G$ is meant in the sense of de la Harpe \cite{dlH2}: every continuous action of $G$ on a compact space admits an invariant regular Borel probability measure. Equivalently, the space $\RUCB(G)$ of bounded right uniformly continuous functions on $G$ admits a left-invariant mean. The right uniform continuity of a function $f$ is meant in the sense of Bourbaki (\cite{bourbaki}, Ch. III, p. 20):  $\abs{f(x)-f(vx)}<\e$ for every $\e>0$ when $v$ belongs to a sufficiently small neighbourhood of identity, $V$. In other words, the map sending each $g\in G$ to the function $\,^gf$, $\,^gf(x)=f(g^{-1}x)$, is continuous in the supremum norm.

The question was advertised by its authors since at least 2000, as a possible tool to prove the existence of invariant vacuum states in gauge field theories. More generally, it can be asked about the group of vertical automorphisms of a smooth principal $K$-bundle for any compact connected Lie group $K$. It can be shown that the $C^{\infty}$ case is equivalent to the amenability of the $C^k$ class groups for all $k\geq 0$, or of the class $H^k$ groups for all $k>\dim X/2$, where $X$ is the base. To our knowledge, the problem remains open for all values of $k$.

An exception is the case of continuous paths and loops (that is, $k=0$ and $X=[0,1]$ or $\s^1$). Here a positive answer follows from the work of Marie-Paule Malliavin e Paul Malliavin \cite{MM}. (See Sect. \ref{s:mm} below.)

In this Note we study groups of finite energy paths and loops (those of Sobolev class $H^1=W^1_2$, strictly intermediate between $C^0$ and $C^1$).
Here is the main result.

\vskip .2cm
\begin{theorem}
Let $K$ be a compact connected Lie group. The topological groups of finite energy paths $H^1(\I,K)$, loops $H^1(\s^1,K)$, based paths $H_0^1(\I,K)$, and based loops $H_0^1(\s^1,K)$, as well as their central extensions, admit a left-invariant mean on the space of left uniformly continuous bounded functions.
\label{th:main}
\end{theorem}
\vskip .2cm

A function $f$ is {\em left uniformly continuous} if for every $\e>0$, we have $\abs{f(x)-f(xv)}<\e$ when $v\in V$ and $V$ is a sufficiently small neighbourhood of identity. In the case of a metrizable group, it is the same as uniform continuity with regard to a compatible left-invariant metric.

We adopt a suggestion by Martin Schneider to call a topological group $G$ admitting a left-invariant mean on the space $\LUCB(G)$ {\em skew-amenable}. 
We still do not know whether the groups of $H^1$ paths and loops are amenable. However, it appears that the skew-amenability of a topological group allows to construct invariant states on algebras of operators with more ease than amenability.

Given a strongly continuous unitary representation $\pi$ of $G$ in a Hilbert space $\H$, a state $\phi$ on the algebra $B(\H)$ (that is, a positive linear functional with $\phi(1)=1$) is invariant if $\phi(\pi_{g^\ast}T\pi_g)=\phi(T)$ for every bounded operator $T$ and each $g\in G$. A unitary representation admitting an invariant state is called {\em amenable} in the sense of Bekka \cite{Bek1}.

\vskip .2cm
\begin{proposition}[Giordano and Pestov \cite{GP07}, Prop. 4.5]
Every strongly continuous unitary representation of a skew-amenable topological group $G$ is amenable in the sense of Bekka.
\label{p:gpskew}
\end{proposition}
\vskip .2cm

For a locally compact group $G$ skew-amenability is equivalent to amenability: since the left regular representation of $G$ is amenable by the above, the group is amenable by Thm. 2.2 in \cite{Bek1}. (For a stronger result, see \cite{greenleaf}, Thm. 2.2.1.) At the same time, the unitary group $U(\ell^2)$ with the strong operator topology is amenable \cite{dlH2} but not skew-amenable, by Prop. \ref{p:gpskew} above and the proof of Prop. 2 in \cite{dlH1} or Ex. 3.6.3 in \cite{P06}. We do not know if every skew-amenable group is amenable.
 
\vskip .2cm
\begin{corollary}
Let $K$ be a compact connected Lie group. The groups of finite energy paths and loops with values in $K$, as well as their central extensions, admit an invariant state for every strongly continuous unitary representation.
\end{corollary}
\vskip .2cm

Note that in a toy example discussed in \cite{CG}, Sect. 4, the group of gauge transformations ($C^0((0,1),SU(n))$ with the relative weak topology) is not only amenable, but skew-amenable as well, because it is a SIN group (the right and left uniform structures coincide). We suggest the following.

\vskip .2cm
\begin{question} Let $K$ be a compact Lie group, and $P$ a principal smooth $G$-fibre bundle. Are the groups of vertical automorphisms of $P$ of classes $C^k$, $H^k$, $k\geq 0$, and $C^{\infty}$ with the corresponding topologies skew-amenable?
\end{question}
\vskip .2cm

\subsection*{Acknowledgements}
The results were announced at the Second Brazilian Workshop in Geometry of Banach Spaces, Ubatuba, SP, Brazil, August 2018. The author is thankful to the Workshop Organizers for an invitation to the extremely pleasant event, and to F. Martin Schneider for stimulating conversations. Support from CNPq (bolsa Pesquisador Visitante, processo 310012/2016) and CAPES (bolsa Professor Visitante Estrangeiro S\^enior, processo 88881.117018/2016-01) is gratefully acknowledged.

\section{Continuous paths and loops: result of Malliavin and Malliavin\label{s:mm}}

According to \cite{MM}, on the group of $C^0$-loops taking values in a compact connected Lie group $K$, integration with regard to the nonlinear Wiener measure $w_t$ of parameter $t$ is asymptotically left-invariant in the limit $t\to\infty$ under the action of the subgroup of $C^1$-loops. In other words, if $f$ is a bounded Borel function on the group $C^0(\s^1,K)$, then for every $C^1$ loop $g$,
\[\int f(g^{-1}x)\,d w_t(x) - \int f(x)\,dw_t(x) \to 0\mbox{ as }t\to\infty.\]
(Apply Theorem 1.1, {\em loco citato}, with $p=2$, and the Cauchy inequality.)
Let now $\mathcal U$ be an untrafilter on $\R$ containing every interval $[t,+\infty)$. It follows that the ultralimit
\[\phi(f) = \lim_{t\to{\mathcal U}}\int f(x)\,dw_t(x)\]
is a $C^1$-left invariant mean on the space of bounded Borel functons on the group of $C^0$-loops. The definition of right uniform continuity easily implies that the restriction of $\phi$ to the bounded right uniformly continuous functions is already invariant under left translations by all $C^0$-loops. So the group $C^0(\s^1,K)$ is amenable. 

A similar argument holds for the path group, but is unnecessary, because the group $C^0(\I,K)$ is isomorphic to a topological quotient group of the loop group. We have:

\vskip .2cm
\begin{corollary}
The groups of continuous paths and loops with values in a compact connected Lie group, with the $C^0$-topology, are amenable.
\end{corollary}

\section{Finite energy paths and loops}

Let $G$ be a Lie group, and $t\in \I$ a point of smoothness of a function $f\colon\I\to G$. The {\em right logarithmic derivative} of $f$ at $t$ is an element of the Lie algebra $\mathfrak g$ of $G$:
\begin{equation}
\partial^{log}f(t) = f^{\prime}(t)\cdot f(t)^{-1}.
\label{eq:rightlog}
\end{equation}
Here we use a simplified notation for the image of $f^{\prime}(t)\in T_{f(t)}G$ in $T_e(G)={\mathfrak g}$ under the right translation by $f(t)^{-1}$. 
If $G$ is a matrix group, it becomes a genuine product of matrices.

If $t$ is a common point of smoothness for two maps $f,g\colon\I\to G$, a direct calculation verifies the cocycle condition:
\begin{equation}
\partial^{log}(fg)(t)=\partial^{log}f(t) + \Ad_{f(t)}\partial^{log}g(t).
\label{eq:cocycle}
\end{equation}
The formula holds in a very general context of infinite-dimensional Lie groups $G$, see \cite{omori}, p. 250. For matrix groups, the adjoint representation is the conjugation: $\Ad_gy=gyg^{-1}$. 


An absolutely continuous mapping $f\colon \I\to G$ has the right logarithmic derivative defined almost everywhere. Such a mapping has {\em finite energy} if
\[\int_0^1 \norm{\partial^{log}f(t)}^2dt<\infty.\]


Now assume $G$ to be a compact connected Lie group, realized as a Lie subgroup of $SO(n)$ inside of $M_n(\R)$ with the Hilbert--Schmidt norm.
Since the right hand side in the formula (\ref{eq:rightlog}) becomes a matrix multiplication, and the matrix $f(t)^{-1}$ is orthogonal, it preserves the norm. It follows that for $G$ compact, the finite energy paths are just the mappings $f\colon\I\to G$ of Sobolev class $H^1=W^{(1)}_2$. 

The space $H^1(\I,M_n(\R))\cong H^1[0,1]\otimes M_n(\R)\cong M_n(H^1[0,1])$ with the Sobolev norm 
\begin{equation}
\norm{f} = \left(\norm{f}_2^2+ \norm{f^{\prime}}^2_2 \right)^{1/2}
\label{eq:sobolev}
\end{equation}
forms a unital Banach algebra under the pointwise multiplication (because so does $H^1[0,1]$ \cite{palais}, Corol. 9.7, p. 30). Consequently, the collection $H^1(\I,G)$ forms a topological subgroup of the group of invertible elements $H^1(\I,M_n(\R))^{\times}$. (It is even a Banach--Lie group, see a more general result in \cite{AHKMT}, Sect. 1.8, but we will not use this.) In the same way,  finite energy loops form a topological group, $H^1(\s^1,G)$, isomorphic to a closed cocompact subgroup of $H^1(\I,G)$.

On the subgroup $H^1_0(\I,G)$ of based paths ($f(0)=e_G$), the Maurer-Cartan cocycle map $\partial^{log}\colon H_0^1(\I,G)\to L^2 (\I,{\mathfrak g})$ is clearly injective. It is also surjective, with the inverse given by the product integral \cite{DF,omori}. 
If $f$ is a step function taking constant values in $\mathfrak g$ on each interval $[t_i,t_{i+1})$, $i=0,1,\ldots,n-1$, and if $t\in [t_j,t_{j+1})$, the product integral of$f$ between $0$ and $t$ is defined by
\[\prod_0^t\exp f(s)\,ds =\exp_G (t-t_j)f(t_j)\exp_G(t_j-t_{j-1})f(t_{j-1})\ldots\exp_G t_1f(0).\]
Afterwards, the product integral is extended over all $L^1$-functions by continuity. In particular, it is well defined on the $L^2$ functions. The function $\I\ni t\mapsto \prod_0^t\exp f(s)\,ds\in G$ is absolutely continuous and satisfies 
\[\partial_t^{log}\prod_0^t\exp f(s)\,ds = f(t),\mbox{ $t$-a.e.}\]
(See \cite{DF}, Thm. 1.2, and Sect. 1.8, p. 55, eq. (8.6).)

Notice that for $f\in H^1_0(\I,G)$, $\norm{f^\prime}_2=\norm{\partial^{log}f}_2$. It follows that on the group $H^1_0(\I,G)$, the
metric induced by the Sobolev norm in Eq. (\ref{eq:sobolev}) is equivalent to the metric induced from the space $L^2((0,1),{\mathfrak g})$ through the cocycle map:
\[d(f,g)=\norm{\partial^{log}f-\partial^{log}g}.\]
This metric $d$ is left invariant, because the $L^2$ norm is invariant under translations and rotations.

One can now identify $H^1_0(\I,G)$ as a topological group with the space $L^2((0,1),\mathfrak g)$ equipped with the multiplication operation
\begin{equation}
f\ast g = f + \Ad_{\Pi\exp f}g.
\label{eq:*}
\end{equation}
Here $\Pi\exp f$ is the short for the function $t\mapsto \prod_{0}^t \exp f(s)\,ds$, and the adjoint representation $\Ad$ is applied pointwise. In our case of $G$ being a matrix group, $\Ad_gx=gxg^{-1}$.

Under the topological group isomorphism
\begin{equation}
\partial^{log}\colon H^1_0(\I,G)\to \left(L^2((0,1),{\mathfrak g}),\ast\right),
\label{eq:isomorphism}
\end{equation}
the additive uniform structure on $L^2((0,1),\mathfrak g)$ corresponds to the left uniform structure on $H^1_0(\I,G)$.
It follows that the group $H^1_0(\I,G)$ is separable and Weil-complete, that is, complete in the left (or, equivalently, right) uniformity. In particular, it is a Polish group. Since $H^1(\I,G)$ is an extension of the based path group by a compact subgroup, and the group $H^1(\I,G)$ contains $H^1(\s^1,G)$ and $H_0^1(\s^1,G)$ as closed subgroups, they are separable and Weil-complete too. 
At the same time, these are not SIN groups.
\vskip .2cm
\begin{proposition} Let $G$ be a non-abelian connected compact Lie group. On each of the groups $H^1_0(\I,G)$, $H^1(\I,G)$, $H^1_0(\s^1,G)$ and $H^1(\s^1,G)$
the left and the right uniform structures are different.
\end{proposition}

\begin{pf} 
Enough to prove the result for the former group (which is obtained from every other group by taking subgroups and/or quotients).
Fix any element $y\in{\mathfrak g}\setminus Z({\mathfrak g})$. Identify $H^1_0(\I,G)$ as a topological group with $(L^2((0,1),\mathfrak g),\ast)$ as in Eqs. (\ref{eq:*}-\ref{eq:isomorphism}). Let $\pi$ be the orthogonal projection of $L^2((0,1),\mathfrak g)$ onto the one-dimensional subspace spanned by the constant function $y$. This $\pi$ is (left) uniformly continuous. Given $\e>0$, there is $f\in L^2((0,1),\mathfrak g)$ with $\norm{f}_2<\e$ and such that $\Ad_{\prod_0^t\exp f(s)\,ds}y\neq \pm y$ for almost all $t$. Since the adjoint action is by the orthogonal transformations, it follows that $\left\Vert\pi\left(\Ad_{\prod\exp f}y - y\right)\right\Vert>0$. 
Now for a sufficiently large $R>0$, 
\begin{align*}
\norm{\pi(f \ast Ry)-Ry}_2& \geq 
R\left\Vert{\pi\left(\Ad_{\prod\exp f}y - y\right)}\right\Vert_2-\e
\end{align*} will be as large as desired. So, $\pi$ is not right uniformly continuous.  \qed
\end{pf}

\section{Skew-amenability of the group of based paths} 

\begin{lemma}
Let $s_n$ be a sequence of positive reals with $s_n=o(n^{-1/2})$. Let $B_n$ denote the unit ball in the $n$-dimensional Euclidean space, and $\lambda_n$ the Lebesgue measure. Then $\lambda^n((B_n+s_ne_1)\setminus B_n)/\lambda^n(B_n)\to 0$ as $n\to\infty$.
\label{l:osqrtn}
\end{lemma}

\begin{pf} 
For $i=n,n-1$ denote $\tilde\lambda_i$ the Lebesgue measure on the ball $B_i$ normalized to one.
Since $s_n^2=o(1/n)$, it is possible to choose a sequence of positive reals $\alpha_n$ with
\[\alpha_n=o(1/n)\cap \omega(s_n^2).\]
The normalized $(n-1)$-volume of the ball $(1-\alpha_n)B_{n-1}$ equals $(1-\alpha_n)^{n-1}=1-o(1)$. For every $x\in e_1^{\perp}$ with $\norm x\leq 1-\alpha_n$, the height $h_n(x)$ of the section of the ball $B_n$ in the direction $e_1$ over the point $x$ satisfies $h_n(x)^2+(1-\alpha_n)^2=1$, thus $h_n(x)^2 =2\alpha_n-\alpha_n^2$ and $h_n(x)=\Omega(\alpha_n^{1/2})$. In particular, $h_n(x)=\omega(s_n)$, uniformly in $x$. Using the Fubini theorem, 
\begin{align*}
\frac{\lambda^n((B_n+s_ne_1)\setminus B_n)}{\lambda^n(B_n)} & = \tilde\lambda_n((B_n+s_ne_1)\setminus B_n) \\
& \leq 
\int_{B_{n-1}\setminus (1-\alpha_n)B_{n-1}}d\tilde\lambda^{n-1}(x) +
\int_{(1-\alpha_n)B_{n-1}} \frac{s_n}{h_n(x)}\,d\tilde\lambda^{n-1}(x) \\
&= o(1) + (1-o(1))o(1) \\
&=o(1).
\end{align*}
\qed
\end{pf}

We will freely use the isomorphism in Eq. (\ref{eq:isomorphism}).
For every $N$, denote $P_N$ the uniform partition of the unit interval into subintervals of length $1/N$. Let $V_N$ be the space of all functions in $L^2((0,1),{\mathfrak g})$ constant a.e. on elements of $P_N$. The union
\[V_{\infty}=\bigcup_{N=1}^{\infty} V_N\]
forms a vector space which is dense in $L^2((0,1),{\mathfrak g})$.

Let $R_N\uparrow\infty$ be a sequence of positive real numbers. Denote $B_N$ the unit ball in the $Nd$-dimensional Euclidean space $V_N$, where $d=\dim{\mathfrak g}$. Let $\nu_N$ be the Lebesgue measure on the ball $R_NB_N$ normalized to one. 

Let $\mathcal U$ be a non-principal ultrafilter on the natural numbers. Define a mean $\phi$ on the space of all bounded (left) uniformly continuous functions $F$ on $L^2((0,1),{\mathfrak g})$ by
 \[\phi(F) = \lim_{N\to{\mathcal U}}\int F(x)\,d\nu_N(x).\]
We will show that, under suitable assumptions on the rate of growth of $R_N$, the mean $\phi$ is left-invariant with respect to the group operation $\ast$
from Eq. (\ref{eq:*}). It is enough to show separately that for every $g\in L^2((0,1),{\mathfrak g})$ we have $\int (F(x+g)-F(x))\,d\nu_N(x)\to 0$, and for every $r\in H^1_0(\I,G)$, $\int (F(\Ad_rx)-F(x))\,d\nu_N(x)\to 0$, as $N\to\infty$.
\vskip .2cm

\begin{lemma}
Suppose 
\[R_N=\omega(N^{1/2}).\]
Let $F$ be a bounded Borel measurable function on $L^2((0,1),{\mathfrak g})$, and let $(g_N)$ be a norm-bounded sequence, where $g_N\in V_N$ for each $N$.
Then 
\[\int \left[F(g_N+ f)-F(f)\right]\,d\nu_N(f)\to 0,\]
and the rate of convergence only depends on $\sup_{N}\norm{g_N}$ and $\norm{F}_{\infty}$. 
\end{lemma}

\begin{pf} 
If every element of the sequence $R_NB_N$ is seen as a F\o lner set in the additive group of the finite-dimensional Euclidean space $V_N$, then the F\o lner constants only depend on $N$ as long as the elements $g_N\in V_{N}$ have norm $\leq M$. In other words, the measure of the symmetric difference between the unit ball $B_N$ and its shift by a vector of length $\leq M$ in whatever direction converges to zero as $N\to\infty$. This is assured by Lemma \ref{l:osqrtn}:
\begin{align*}
\int_{L^2} \left[F(g_N +f)-F(f)\right]\,d\nu_N(f) & =
\int_{R_NB_N}F(g_N+ f)\,d\nu_N(f)-\int_{R_NB_N}F(f)\,d\nu_N(f) \\
&= \int_{(R_NB_N+g_N)\cap R_NB_N}F(f)\,d\nu_N(f)+
\int_{R_NB_N\setminus R_NB_N-g_N}F(f)\,d\nu_N(f)\\
&\phantom{xxx}
 -\int_{(R_NB_N+g_N)\cap R_NB_N}F(f)\,d\nu_N(f) - \int_{R_NB_N\setminus (R_NB_N+g_N)}F(f)\,d\nu_N(f)\\
&\leq 2\norm{F}_{\infty}\nu_N(R_NB_N\setminus g_N+R_NB_N) \\
&= 2\norm{F}_{\infty}\tilde\lambda_n (B_N\setminus \norm{g_N}R_N^{-1}e_1+B_N) \\
&\leq 2\norm{F}_{\infty}\tilde\lambda_n (B_N\setminus MR_N^{-1}e_1+B_N) \\
&\to 0,
\end{align*}
because $MR^{-1}_N=o(N^{-1/2})$ by assumption. \qed
\end{pf}

Recall that for an element $f\in L^2((0,1),{\mathfrak g})$ and a function $r\colon [0,1]\to G$, we define $\Ad_rf$ pointwise a.e., that is, $(\Ad_rf)(x)=\Ad_{r(x)}f(x)$. Since we realize $\mathfrak g$ and $G$ inside of $M_n(\R)$, the adjoint representation is just the conjugation: $\Ad_g(x)=gxg^{-1}$. Denoting the uniform operator norm on matrices $\norm{\cdot}_u$, and the Hilbert--Schmidt norm by $\norm{\cdot}_2$, we have, for two unitary operators $u,v$ and a matrix $A\in M_n$,
\begin{align*}
\norm{\Ad_uA-\Ad_vA}_2 &= \norm{uAu^{-1}-vAv^{-1}}_2 \\
& \leq \norm{uAu^{-1}-uAv^{-1}}_2 + \norm{uAv^{-1}-vAv^{-1}}_2 \\
&\leq 2\norm{u-v}_u\norm{A}_2.
\end{align*}
\vskip .2cm

\begin{lemma}
Suppose
\[R_N= o(N/\log N).\]
Then for every $H^1$-function $r\colon [0,1]\to K\subseteq M_n$ and every bounded uniformly continuous function $F$ on $L^2((0,1),{\mathfrak g})$,
\[\int_{L^2} \left[ F(\Ad_rf) - F(f)\right] d\nu_N(f)\to 0\mbox{ as }N\to\infty,\]
and the convergence is uniform in $r$ over every set bounded in the norm $\norm{\partial^{log}r}_2$, and also over any uniformly equicontinuous and uniformly bounded family of functions $F$.
\end{lemma} 

\begin{pf}
Every $H^1$-function $r\colon [0,1]\to G\subseteq M_n$ is H\"older continuous with exponent $1/2$, more exactly, for $0\leq s \leq t\leq 1$,
\begin{align*}
\norm{r(s)-r(t)}_{u} &\leq C\left(\int_s^t \norm{r^{\prime}(x)}_2^2dx\right)^{1/2}\abs{s-t}^{1/2},\\
&=C\left(\int_s^t \norm{r^{\prime}(x)r(x)^{-1}}_2^2dx\right)^{1/2}\abs{s-t}^{1/2}.
\end{align*}
For every $i=0,1,\ldots, N-1$, denote
\[c_i = \left(\int_{i/N}^{(i+1)/N} \norm{r^{\prime}(x)}_2^2dx\right)^{1/2}.\]
Thus, $\sum_{i=0}^{N-1}c_i^2 = \norm{\partial^{\log}r}^2$. Given $t\in [i/N,(i+1)/N)$,
\begin{align*}
\norm{r(i/N)-r(t)}_{u} &\leq Cc_i N^{-1/2}.
\end{align*}
Approximate $r$ with a step function $\rho$ taking the value $r(i/N)$ on the $i$-th interval of the partition. For $t$ as above,
\[\norm{r(t)-\rho(t)}_{u}\leq Cc_i N^{-1/2}.\]
Let $f\in B_N$, the unit ball in the Euclidean space $V_N$. Denote $f_i$ the constant value taken by $f$ on the interval $[i/N,(i+1)/N)$. We have:
\begin{align*}
\norm{ \Ad_rf-\Ad_{\rho} f}_2^2 &= \int_0^1 \norm{(\Ad_{r(x)}-\Ad_{\rho(x)})f(x)}_2^2d x\\
& \leq \sum_{i=0}^{N-1} \int_{i/N}^{(i+1)/N} 4\norm{r(x)-\rho(x)}_u^2\norm{f_i}_2^2d x\\
&\leq 4\sum_{i=0}^{N-1} \int_{i/N}^{(i+1)/N}C^2c_i^2 N^{-1} \norm{f_i}_2^2dx\\
&= C^{\prime}N^{-1}\sum_{i=0}^{N-1}c_i^2  \int_{i/N}^{(i+1)/N} \norm{f_i}_2^2dx\\
&= C^{\prime}N^{-1}\sum_{i=0}^{N-1}c_i^2  N^{-1}\norm{f_i}_2^2\\
&\leq C^{\prime}N^{-1}\norm{\partial^{\log}r}^2
\max_{i=0}^{N-1}\left(N^{-1}\norm{f_i}_2^2 \right).
\end{align*}
Thus,
\begin{equation}
\norm{\Ad_rf-\Ad_{\rho} f}_2\leq C^{\prime\prime}N^{-1/2}\norm{\partial^{\log}r}
\max_{i=0}^{N-1}\left(N^{-1/2}\norm{f_i}_2 \right).
\label{eq:thus}
\end{equation}
For every $i$, the function 
\[B_N\ni f\mapsto N^{-1/2}f_i\in{\mathfrak g}\subseteq M_n\]
is $1$-Lipschitz. Fix an orthonormal basis $e_1,\ldots,e_d$ in $\mathfrak g$. The function $B_N\ni f\mapsto \langle N^{-1/2}f_i,e_j\rangle\in\R$ is $1$-Lipschitz, with zero as the median value. The ball $B_N$ has dimension $Nd$. By the L\'evy concentration inequality for the Euclidean ball (\cite{ledoux}, Prop. 2.9, p. 30), for every $\e>0$
\[\nu\left\{f\in B_N\colon \left\vert\langle N^{-1/2}f_i,e_j\rangle\right\vert>\e\right\}< 2\exp(-c\e^2Nd),\]
where $\nu$ is the Lebesgue measure on $B_N$ normalized to one.
Consequently,
\begin{align*}
\nu\{f\in B_N\colon N^{-1/2}\norm{f_i}_2>\e\} & \leq \nu\left\{x\in B_N\colon \exists j=1,\ldots,d,~~\left\vert\langle N^{-1/2}f_i,e_j\rangle\right\vert>\e\right\} \\
&< 2d\exp(-c\e^2Nd).
\end{align*}
Substituting $\e=\log N/\sqrt N$, 
\begin{align*}
\nu\{x\in B_N\colon N^{-1/2}\norm{f_i}_2>\log N\cdot N^{-1/2}\}&< 2d\exp(-cd\log^2N),
\end{align*}
and
\begin{align*}
\nu\left\{x\in B_N\colon\max_{i=0}^{N-1}\left(N^{-1/2}\norm{f_i}_2 \right)>\log N\cdot N^{-1/2}\right\}& < 2Nd \exp(-cd\log^2N) \\
& = 2Nd N^{-cd\log N} \\
&= 2d N^{1-cd\log N} \\
& = N^{-\omega(1)}.
\end{align*}
Denote
\[A_N = \{f\in B_N\colon \forall i,~~N^{-1/2}\norm{f_i}\leq \log N\cdot N^{-1/2}\}\]
and $A_N^c=B_N\setminus A_N$.

Let $F$ be a uniformly continuous bounded function on $L^2((0,1),{\mathfrak g})$, and denote
\[\e_F(\delta)=\sup\{\abs{F(f)-F(h)}\colon \norm{f-h}_2<\delta\}.\]
We have $\e_F(\delta)\to 0$ when $\delta\to 0$. 

Note that $\Ad_{\rho}$ leaves $V_N$ invariant, and since $\Ad_\rho$ is an orthogonal operator, it leaves the measure $\nu$ invariant as well.  We have 
\begin{align*}
\int_{L^2} \left[F(\Ad_rf) - F(f)\right] d\nu_N(f) &=\int_{R_NB_N} \left[ F(\Ad_rf) - F(\Ad_{\rho}f)\right] d\nu_N(f) +
\int_{R_NB_N} \left[ F(\Ad_{\rho}f) - F(f)\right] d\nu_N(f), \\
\end{align*}
where the second integral vanishes.
Let us estimate the first integral separately over $R_NA_N$ and $R_NA^c_N$. For the latter,
\begin{align*}
\left\vert \int_{R_NA^c_N}\left[ F(\Ad_rf) - F(f)\right] d\nu_N(f) \right\vert
&\leq
\int_{R_NA^c_N} \left\vert F(\Ad_rf) - F(f)\right\vert d\nu_N(f)\\
 &
\leq 2\norm{F}_{\infty}\nu_N(R_NA^c_N) \\
&= 2\norm{F}_{\infty}\tilde\lambda_N (A^c_N) \\
& < 4d\norm{F}_{\infty} N^{1-cd\log N}  \\
& \to 0,
\end{align*}
and the convergence is obviously uniform in $r$. For $R_NA_N$, we have, using Eq. (\ref{eq:thus}),

\begin{align*}
\sup_{f\in R_NA_N}\norm{\Ad_rf - \Ad_{\rho}f}_2 & \leq 
\sup_{f\in R_NA_N}
N^{-1/2}\norm{\partial^{\log}r}_2\max_{i=0}^{N-1}\left(N^{-1/2}\norm{f_i}_2 \right)
\\
&\leq N^{-1/2}\norm{\partial^{\log}r}_2\frac{\log N}{N^{1/2}}R_N \\
& = \norm{\partial^{\log}r}_2\frac{\log N}{N} R_N \\
&\to 0,
\end{align*}
because of the hypothesis $R_N=o(N/\log N)$. 
 Now,
\begin{align*}
\left\vert\int_{R_NA_N} \left[ F(\Ad_rf) - F(\Ad_{\rho}f)\right] d\nu_N(f) \right\vert &\leq
\int_{R_NA_N} \left\vert F(\Ad_rf) - F(\Ad_{\rho}f)\right\vert d\nu_N(f) \\
& \leq\e_F\left(\sup_{f\in R_NA_N}\norm{\Ad_rf - \Ad_{\rho}f} \right) \\
&\leq \e_F\left(\norm{\partial^{\log}r}_2\frac{\log N}{N} R_N\right)\\
&\to 0.
\end{align*}  The convergence is uniform in $r$ on any set where the values of $\norm{\partial^{\log}r}_2$ are bounded. It is also uniform for any uniformly bounded family of functions $F$ having a common modulus of uniform continuity. \qed
\end{pf}

\vskip .2cm
\begin{lemma}
Suppose $R_N=o(N/\log N)\cap\omega(N^{1/2})$.
Let $D$ be a compact subset of the group $H^1_0((0,1),G)$. For every bounded left uniformly continuous function $F$ on the group we have, uniformly in $g\in D$,
\begin{equation}
\int \left[F(g\ast f) - F(f)\right]d\nu_N(f) \to 0.
\label{eq:therate}
\end{equation}
The convergence is uniform over any uniformly equicontinuous and uniformly bounded family of functions $F$.
\label{l:skewH1}
\end{lemma}

\begin{pf} 
Denote $M=\max_{f\in D}\norm{f}+1$. Given $\e>0$, find
$\delta>0$ such that whenever $\norm{f-h}_2<\delta$, $\abs{F(f)-F(h)}<\e$.
Now find $N$ so large that, for each $n\geq N$,
\vskip .2cm

\begin{itemize}
\item $D$ is contained in the $\delta$-neighbourhood of $V_n$,
\item $\int\left[F(g+f)-F(f) \right]d\nu_n(f)<\e$ whenever $g^\prime\in V_N$ and $\norm{g^\prime}\leq M$,
\item for each $g\in D$, $\int\left[F(\Ad_{\Pi\exp g}f)-F(f) \right]d\nu_n(f)<\e$.
\end{itemize}
\vskip .2cm

Given $g\in D$ and $n\geq N$, find $g^\prime\in V_n$ with $\norm{g-g^\prime}<\delta$. Now 
\begin{align*}
\left\vert\int \left[ F(g\ast f)-F(f)\right]\,d\nu_n\right\vert &=
\left\vert\int \left[ F(g+\Ad_{\Pi \exp g}f)-F(f)\right]\,d\nu_n \right\vert\\
&\leq\left\vert\int \left[ F(g+\Ad_{\Pi \exp g}f)-F(g^\prime+\Ad_{\Pi\exp  g}f)\right]\,d\nu_n \right\vert + \\ &
\left\vert\int \left[ F(g^\prime+\Ad_{\Pi\exp  g}f)-F(\Ad_{\Pi\exp  g}f)\right]\,d\nu_n \right\vert +\\ &
\left\vert\int \left[ F(\Ad_{\Pi \exp g}f)-F(f)\right]\,d\nu_n \right\vert \\
&<3\e.
\end{align*}
\qed
\end{pf}

As an immediate consequence, we obtain the skew-amenability for the group $H^1_0(\I,K)$ of based paths.

\section{Co-compact normal subgroups}
In this Section we will establish the skew-amenability of the group of based loops, as a direct consequence of the following result.
\vskip .2cm

\begin{theorem}
A normal co-compact subgroup of a skew-amenable metrizable group is skew-amenable.
\end{theorem}
\vskip .2cm

(The result is valid without an assumption of metrizability, but the proof is more involved.)

For a topological group $G$, denote ${\mathcal L}(G)$ the Gelfand space of the $C^\ast$-algebra $\LUCB(G)$. This is the maximal (Samuel) compactification (\cite{samuel}, Ch. III) of the left uniform space of $G$. 
The left action of $G$ on the algebra, $f\mapsto {\,}^gf$, is by $C^\ast$-algebra automorphisms, and so it determines an action of $G$ on the left by homeomorphisms on the space ${\mathcal L}(G)$, although this action will in general be discontinuous. By the Riesz representation theorem, $G$ is skew-amenable if and only if ${\mathcal L}(G)$ supports a left-invariant regular Borel probability measure. 

It is instructive to compare ${\mathcal L}(G)$ with the Samuel compactification of the right uniform space of $G$, known as the {\em greatest ambit}, ${\mathcal S}(G)$ (see e.g. \cite{devries,P06}). In this case, the left action of $G$ on ${\mathcal S}(G)$ is continuous, and $G$ is amenable if and only if ${\mathcal S}(G)$ admits a left-invariant regular Borel probability measure. It is easy to see, using the inversion map $g\mapsto g^{-1}$, that the {\em right} action on the ``skew ambit'' ${\mathcal L}(G)$ is continuous, and $G$ is amenable if and only if ${\mathcal L}(G)$ admits a {\em right} invariant probability measure. Equivalently, $G$ is skew-amenable if the greatest ambit admits a {\em right} invariant probability measure.

Now let $H$ be a normal subgroup of a metrizable topological group $G$, such that the quotient group $K=G/H$ is compact.
Fix a bounded left-invariant compatible metric $d$ on $G$, and denote $\pi\colon G\to G/H$ the quotient homomorphism. 
It extends to a continuous $G$-equivariant map $\bar\pi\colon {\mathcal L}(G)\to G/H$. The formula $\bar d(xH,yH)=\inf\{d(xh_1,yh_2(\colon h_1,h_2\in H\}$ defines a left-invariant compatible metric on $G/H$ (this is where we are using the normality of $H$).
Let $\mu$ be a left-invariant regular Borel probability measure on ${\mathcal L}(G)$.
For $\delta>0$, cover $K$ with finitely many open $\delta$-balls. They are all left translations of each other. The same is true of their inverse images under $\bar\pi$, all of which consequently have the same strictly positive measure. The open set $V_{\delta}=\bar\pi^{-1}B^{\bar d}_\delta(e)$ satisfies $HV_{\delta}=V_{\delta}$, so is invariant under the left translations by elements of $H$. Let $\mu_{\delta}$ be the restriction of the measure $\mu$ to $V_{\delta}$, normalized to one. The measure $\mu_{\delta}$ is a regular Borel probability measure on ${\mathcal L}(G)$, invariant under left translations by elements of $H$.

Every bounded uniformly continuous function on a metric space is uniformly approximated with Lipschitz functions (see e.g. \cite{GJ}, Corol. 1). Therefore, it is enough to define a left-invariant mean on the bounded Lipschitz functions on $(H,d)$, or, essentially the same, a left-invariant positive affine functional on the bounded $1$-Lipschitz functions sending $0$ to $0$ and $1$ to $1$. Given a bounded $1$-Lipschitz function $F\colon H\to\R$, extend it to a bounded $1$-Lipschitz function $\tilde F\colon G\to\R$ in a usual way:
\[\tilde F(g) = \inf\{F(h)+d(g,h)\colon h\in H\}.\]
Further, extend $\tilde F$ to a continuous function $\bar F\colon {\mathcal L}(G)\to \R$. Now choose an ultrafilter $\mathcal U$ on $\R$ containing all the intervals $(0,\delta)$, $\delta>0$, and set
\[\phi(F)=\lim_{\delta\to{\mathcal U}}\int_{{\mathcal L}(G)}\bar F(x)\,d\mu_{\delta}(x).\]

1. {\em Left-invariance of $\phi$.} For every $k\in H$ and $g\in G$,
\begin{align*}
{\,}^k\tilde F(g) &= \tilde F(k^{-1}g) \\
&= \inf\{F(h)+d(k^{-1}g,h)\colon h\in H\}\\
&= \inf\{F(k^{-1}h)+d(k^{-1}g,k^{-1}h)\colon h\in H\}\\
&= \inf\{{\,}^kF(h)+d(g,h)\colon h\in H\}\\
&= \widetilde{{\,}^kF}(g).
\end{align*}
By continuity, one must have ${\,}^k\bar F=\overline{{\,}^kF}$, and since the measures $\mu_{\delta}$ are left-invariant, we conclude.
\vskip .2cm

2. {\em Other properties of $\phi$.} Given $x\in V_{\delta}\cap G$, there is $h\in H$ with $d(h,x)<\delta$. It means that the value at $x$ of every function $\tilde F$, where $F$ is $1$-Lipschitz on $H$, is $\delta$-approximated by the value of $F$ at $h$. Consequently, for every constant function $c$ on $H$, $\sup_{x\in V_{\delta}}\left\vert \bar c(x)-c \right\vert\leq \delta$, and since $\mbox{supp}\,\mu_{\delta}\subseteq V_{\delta}$,
\[\int \left\vert \bar c(x)-c \right\vert\,d\mu_{\delta}(x)\leq \delta.\]
 For every positive $1$-Lipschitz function $F\geq 0$, we deduce that $\tilde F\geq -\delta$ on $V_{\delta}$, so
\[\int \bar F(x) \,d\mu_{\delta}(x)\geq -\delta.\]
 For any $1$-Lipschitz functions $F,G$ and $t\in [0,1]$, we have 
\[\sup_{x\in V_{\delta}}\left\vert \widetilde{tF+(1-t)G}(x) - t\tilde F(x)-(1-t)\tilde G(x)\right\vert \leq 3\delta,\]
and 
\[\int \left\vert\widetilde{tF+(1-t)G}(x) - t\tilde F(x)-(1-t)\tilde G(x)\right\vert\,d\mu_{\delta}\leq 3\delta.\]
In the ultralimit, we obtain the desired properties of $\phi$. \qed

\section{Semidirect products with compact groups}

The groups of paths $H^1(\I,G)$ and loops $H^1(\s^1,G)$ decompose into semidirect products of the normal subgroups, $H_0^1(\I,G)$ and $H_0^1(\s^1,G)$ respectively, with the compact subgroup of constant paths ($=$ loops), isomorphic to $G$. Therefore, to establish their skew-amenability, it is enough to prove to following result.
\vskip.2cm

\begin{theorem}
The semidirect product of a compact group and a skew-amenable group is  skew-amenable.
\label{l:semidirect}
\end{theorem}
\vskip .2cm

Recall that a topological group $G$ is an {\em inner} semidirect product of a subgroup $K$ and a normal subgroup $H$ if every element of $G$ is uniquely written as a product of an element of $K$ and an element of $H$, and the resulting map $K\times H\to G$ is a homeomorphism. In this case, 
for $x,y\in K$ and $g,h\in H$,
\[xgyh = xy\cdot y^{-1}gy\cdot h,\]
and the mapping $\varsigma\colon K\to \Aut H$, $\varsigma_yg=y^{-1}gy$, is an antihomomorphism, that is, $\varsigma_{yz}g=z^{-1}y^{-1}gyz =\varsigma_{z}\varsigma_{y}g$. The group operation becomes
\[(x,g)(y,h)=(xy,\varsigma_yg\cdot h).\]
We will denote this relization of outer semidirect product $K\ltimes^1_\varsigma H$.

Note however that the group $G$ also admits a unique decomposition as $G=H\cdot K$: indeed, $g= xh =xhx^{-1}\cdot x$, and the mapping $(h,x)\mapsto hx = x\cdot x^{-1}hx\in G$ is also a homeomorphism. In this case, 
\[gxhy = g\cdot xhx^{-1}\cdot xy,\]
and the map of $K$ to $\Aut H$ defined by $\tau_xh= xhx^{-1}$ is a homomorphism. The group operation on the product $K\times H$ becomes
\[(x,g)(y,h)=(xy, g\tau_xh).\]
We will denote this realization of semidirect product $K\ltimes^2_\tau H$.
Here is an isomorphism between the two:
\begin{equation}
K\ltimes^1_\tau H\ni (x,h)=xh\mapsto xhx^{-1}\cdot x=(x,\tau_xh)\in K\ltimes^2_\tau H.
\label{eq:isomorphism2semidirect}
\end{equation}
 
\vskip .2cm

\begin{lemma}
Under our assumptions, the compactification ${\mathcal L}(K\ltimes^2_{\tau} G)$ can be canonically identified with the product $K\times {\mathcal L}(G)$: the unique continuous map ${\mathcal L}(K\ltimes^2_{\tau} G)\to K\times {\mathcal L}(G)$ extending the identity map on $K\ltimes^2_{\tau} G$ is a homeomorphism.
\end{lemma} 

\begin{pf}
Let $f\colon K\ltimes^2 G\to\R$ be a bounded uniformly continuous function. 
For every $x\in K$, the function 
\[f_{x}\colon G\ni g\mapsto f(x,g)\in\R\]
is bounded left uniformly continuous on $G$, and so extends uniquely to a continuous function $\overline{f_x}\colon {\mathcal L}(G)\to\R$. Coalescing the functions $\overline{f_x}$, $x\in K$, we get a bounded function
\[\bar f\colon K\times {\mathcal L}(G)\ni (x,g)\mapsto \overline{f_x}(g)\in\R,\]
so it only remains to show it is continuous. 

Let $(x_0,\xi_0)\in K\times {\mathcal L}(G)$, and let $\e>0$.
Since the right action by $K$ on ${\mathcal L}(G)$ is continuous, there is a neighbourhood of identity, $V$, in $K$, with the property: for all $x\in V$ and all $(y,h)\in K\ltimes G$, 
\[f(yx,h)\overset\e\approx f(y,h)\]
(here $(y,h)(x,e)=(yx,h)$). Let $\zeta\in {\mathcal L}(G)$ be arbitrary. There exist a net of elements of $G$ with $g_{\alpha}\to\zeta$. For every $x\in K$, $\bar f(x,g_{\alpha})=\overline{f_x}(g_{\alpha})\to \overline{f_x}(\zeta)=\bar f(x,\zeta)$. We conclude: for all $x\in V$ and all $(y,\zeta)\in K\ltimes {\mathcal L}(G)$,
\begin{equation}
\bar f(yx,\zeta)\overset\e\approx \bar f(y,\zeta).
\label{eq:xi}
\end{equation}
Now define an open neighbourhood of $\xi_0$ in ${\mathcal L}(G)$ by
\begin{equation}
U=\left(\overline{f_{x_0}}\right)^{-1}(f(x_0,\xi_0)-\e,f(x_0,\xi_0)+\e).
\label{eq:U}
\end{equation}
The set $x_0V\times U$ is an open neighbourhood of $(x_0,\xi_0)$ in $K\times {\mathcal L}(G)$.
Let $(x,\xi)\in x_0V\times U$. For some $v\in V$, we have $x=x_0v$, and 
\begin{align*}
\left\vert\bar f(x,\xi)-\bar f(x_0,\xi_0)\right\vert &\leq
\left\vert\bar f(x_0v,\xi)-\bar f(x_0,\xi)\right\vert +
\left\vert\bar f(x_0,\xi)-\bar f(x_0,\xi_0)\right\vert \\
&< 2\e,
\end{align*}
by Eqs. (\ref{eq:xi}) and (\ref{eq:U}). 
\qed
\end{pf}

Observe that the compactness of $K$ in the lemma above cannot be dropped: for an infinite discrete group $G$, the ``skew ambit'' ${\mathcal L}(G)$ is just the Stone--\v Cech compactification, $\beta G$, and it is well known that $\beta(\N\times \N)\neq \beta\N\times \beta\N$ (\cite{engelking}, Exercise 3.6.D(b)). In fact, the Stone--\v Cech compactification of the product of a compact group $K$ with a discrete group $D$ need not be equal to the product of $K$ with $\beta D$ either (\cite{engelking}, Exercises 3.6.D(a) and 3.2.H(b)).
\vskip .2cm

Every topological group automorphism $\varsigma_x$, $x\in K$ of $G$ extends to a $C^\ast$-algebra automorphism of $\LUCB(G)$, and further to a homeomorphism of ${\mathcal L}(G)$. We will denote it with the same symbol, $\varsigma_x$. 

\vskip .2cm 
\begin{lemma}
Under our assumptions, the compactification ${\mathcal L}(K\ltimes^1_{\tau} G)$ can be canonically identified with the product $K\times {\mathcal L}(G)$, with the action of the group $K\ltimes^1_{\tau} G$ given by
\begin{equation}
(x,g)(y,\xi)=(xy,\varsigma_y g\cdot \xi),
\label{eq:formulasemidirect}
\end{equation}
for all $x,y\in K$, $g\in G$, and $\xi\in {\mathcal L}(G)$.
\end{lemma} 

\begin{pf}
The group isomorphism in Eq. (\ref{eq:isomorphism2semidirect}) preserves the left uniform structure, sends the subgroup $K$ to itself, and preserves every fiber $\{x\}\times G$, from where the first statement follows. 

To prove the formula, choose a net $g_{\alpha}$ of elements of $G$ converging to $\xi$.  For fixed $x\in K$ and $g\in G$, the left translations of ${\mathcal L}(K\ltimes G)\cong K\times {\mathcal L}(G)$ by elements of the group $K\ltimes G$ are continuous, and so are the transformations $\tau_x$ and the multiplication on the left on ${\mathcal L}(G)$, therefore
 \begin{align*}
(x,g)(y,\xi)& = (x,g)(y,\lim_{\alpha}\, g_{\alpha})\\
&= (x,g)\lim_{\alpha}\,(y, g_{\alpha})\\
&= \lim_{\alpha}\, (x,g)(y,g_{\alpha}) \\
&= \lim_{\alpha}\, (xy,\varsigma_y g\cdot g_{\alpha}) \\
&= (xy,\lim_{\alpha}\, \varsigma_y g \cdot g_{\alpha}) \\
& = (xy,\varsigma_yg \lim_{\alpha}\, g_{\alpha}) \\
&=(xy,\varsigma_y g\xi).
\end{align*}
To get the equalities on lines 2 and 6, we have implicitely used the topological identification ${\mathcal L}(K\ltimes^1_{\tau} G)\cong K\times {\mathcal L}(G)$.
\qed
\end{pf}

Now we are ready to prove Thm. \ref{l:semidirect}. 
Fix a left-invariant regular Borel probability measure $\mu$ on ${\mathcal L}(G)$, and let $\nu$ be the normalized Haar measure on $K$. We claim that the  product probability measure $\nu\otimes\mu$ on $K\times^1_{\varsigma} {\mathcal L}(G)$ is left-invariant. Indeed, for every continuous function $f$ on $K\times {\mathcal L}(G)$ and each $(\kappa^\prime,g )\in K\ltimes^1_{\varsigma} G$, we have, using the Fubini theorem (\cite{fremlin2}, Th. 252B):

\begin{align*}
\int {\,}^{(\kappa^\prime,g )}f(\kappa,s)\,d(\nu\otimes\mu)(\kappa,s) &= \int_Kd\nu(\kappa) 
\int_{{\mathcal L}(G)} f((\kappa^\prime,g )^{-1}\cdot(\kappa,s))\,d\mu(s)\\
&= \int_Kd\nu(\kappa) \int_{{\mathcal L}(G)} 
f((\kappa^{\prime (-1)},\tau_{\kappa^{\prime (-1)}}g^{-1})\cdot(\kappa,s)))\,d\mu(s) \\
&=\int_Kd\nu(\kappa) \int_{{\mathcal L}(G)} 
f\left(\kappa^{\prime (-1)}\kappa, \tau_{\kappa^{\prime (-1)}\kappa}g^{-1}\cdot s \right) \,d\mu(s) \\
&=\int_Kd\nu(\kappa) \int_{{\mathcal L}(G)} 
f\left(\kappa^{\prime (-1)}\kappa, s \right) \,d\mu(s) \\
&=\int_Kd\nu(\kappa)\,\, {\,}^{\kappa^{\prime}}\!\!\left[\kappa\mapsto \int_{{\mathcal L}(G)} f(\kappa,s)\,d\mu(s) \right]\\
&=\int_Kd\nu(\kappa) \int_{{\mathcal L}(G)} f(\kappa,s)\,d\mu(s) \\
&=\int f(\kappa,s)\,d(\nu\otimes\mu)(\kappa,s).
\end{align*}
\qed

\section{Central extensions}

\begin{proposition}
A central extension of a skew-amenable topological group is skew-amenable.
\label{p:central}
\end{proposition}

\begin{pf}
Let $e\to Z\to G\to G/Z\to e$, where $Z$ is central and $G/Z$ is skew-amenable. Select an invariant mean $\phi_Z$ on $\ell^{\infty}(Z)$ and a left-invariant mean, $\phi_{G/Z}$, on $\LUCB(G/Z)$. Given a bounded left-uniformly continuous function $f$ on $G$, define
\[\phi_G(f) = \phi_{G/Z}\left[gZ\mapsto \phi_Z\left[z\mapsto f(zg)\vert_Z \right] \right].
\]
The intermediate function $gZ\mapsto \phi_Z\left[z\mapsto f(zg)\vert_Z \right]$ is well-defined (that is, constant on the $Z$-cosets), because for each $z^\prime\in Z$,
\begin{align*}
\phi_Z\left[z\mapsto f(zgz^\prime)\vert_Z \right] &= \phi_Z\left[z\mapsto f(z^\prime zg)\vert_Z \right]\\
&= \phi_Z\left[z\mapsto {\,}^{z^{\prime{-1}}}f(zg)\vert_Z \right]\\
&= \phi_Z\left[z\mapsto f(zg)\vert_Z\right].
\end{align*}
The intermediate function is also obviously bounded, and left uniformly continuous: if $V$ is such that $g^{-1}h\in V$ implies $f(g)\overset\e\approx f(h)$, then for all $z\in Z$, also $f(zg)\overset\e\approx f(zh)$, and $\phi_Z\left[z\mapsto f(zg) \right]\overset\e\approx \phi_Z\left[z\mapsto f(zh) \right]$. Thus, the value $\phi_G(f)$ above is well-defined, and is clearly a mean on $\LUCB(G)$. Remains to verify the left-invariance:
\begin{align*}
\phi_G(^hf)&
 = \phi_{G/Z}\left[gZ\mapsto \phi_Z\left[z\mapsto {\,}^hf(zg)\vert_Z \right] \right]\\
&= \phi_{G/Z}\left[gZ\mapsto \phi_Z\left[z\mapsto f(h^{-1}zg)\vert_Z \right] \right]\\
& = \phi_{G/Z}\left[gZ\mapsto \phi_Z\left[z\mapsto f(zh^{-1}g)\vert_Z \right] \right]\\
&= \phi_{G/Z}\left(^{hZ}\left[gZ\mapsto \phi_Z\left[z\mapsto f(zg)\vert_Z \right] \right]\right)\\
&= \phi_{G/Z}\left[gZ\mapsto \phi_Z\left[z\mapsto f(zg)\vert_Z \right] \right]\\
&=\phi_G(f).
\end{align*}
\qed
\end{pf}

The argument is of course standard, but uses the centrality of $Z$ in an essential way. We do not know if Prop. \ref{p:central} and Thm. \ref{l:semidirect} can be generalized to show that skew-amenability is closed under extensions. 

In the non-locally compact case, skew-amenability appears to be more fragile than amenability, because of discontinuity of the left action. It could certainly help, to have a general criterion of skew-amenability as powerful as the criterion of amenability for all topological groups established recently by Schneider and Thom (Thm. 3.2 in \cite{ST}). According to a private communication from Martin Schneider, such a result is possible and will be published soon.

\end{document}